\documentclass[11pt, twoside]{article}

\usepackage{extarrows}
\usepackage{enumerate}
\usepackage{amsmath}
\usepackage{amssymb}
\usepackage{amsthm}
\usepackage[all]{xypic}
\usepackage[all,cmtip,poly]{xy}
\usepackage{latexsym}
\usepackage{pstricks,hyperref}
\usepackage{amsbsy}
\usepackage{amscd}
\usepackage{mathrsfs}
\usepackage{txfonts}
\usepackage{amsfonts}
\usepackage{graphicx}
\usepackage[top=1.40in,bottom=1.30in,left=1.20in,right=1.20in]{geometry}
\newtheorem{thm}{Theorem}[section]
\newtheorem{cor}[thm]{Corollary}
\newtheorem{lem}[thm]{Lemma}

\newtheorem{defn}[thm]{Definition}
\newtheorem{rem}[thm]{Remark}

\usepackage{pstricks,color}

\usepackage{indentfirst}
\begin{document}

\begin{center}
{\Large \bf Two-term relative cluster tilting subcategories, $\tau-$tilting modules and silting subcategories\footnotetext{This work was supported by the Hunan Provincial Natural Science Foundation of China (Grants No.2018JJ3205) and the NSF of China (Grants No.11671221)}}

\bigskip

{\large Panyue Zhou and Bin Zhu}
\bigskip

\end{center}

\def\s{\stackrel}
\def\Longrightarrow{{\longrightarrow}}
\def\A{\mathcal{A}}
\def\B{\mathcal{B}}
\def\C{\mathscr{C}}
\def\D{\mathsf{D}}
\def\T{\mathcal{T}}
\def\M{\mathcal{M}}
\def\E{\mathcal{E}}
\def\R{\mathcal{R}}
\def\S{\mathcal{S}}
\def\H{\mathcal{H}}
\def\U{\mathscr{U}}
\def\V{\mathscr{V}}
\def\W{\mathscr{W}}
\def\X{\mathscr{X}}
\def\Y{\mathscr{Y}}
\def\Z{\mathcal {Z}}
\def\I{\mathcal {I}}
\def\RR{\mathcal{R}\ast\mathcal{R}[1]}
\def\Aut{\mbox{Aut}}
\def\coker{\mbox{coker}}
\def\Ker{\mbox{Ker}}
\def\deg{\mbox{deg}}
\def\dim{\mbox{dim}}
\def\Ext{\mbox{Ext}}
\def\Hom{\mbox{Hom}}
\def\DHom{\mbox{DHom}}
\def\Gr{\mbox{Gr}}
\def\id{\mbox{id}}
\def\Im{\mbox{Im}}
\def\ind{\mbox{ind}}
\def\Int{\mbox{Int}}
\def\ggz{\Gamma}
\def\la{\Lambda}
\def\bz{\beta}
\def\az{\alpha}
\def\gz{\gamma}
\def\da{\delta}
\def\fs{{\mathfrak{S}}}
\def\ff{{\mathfrak{F}}}
\def\zz{\zeta}
\def\thz{\theta}
\def\raw{\rightarrow}
\def\ole{\overline}
\def\cat{\C_{F^m}(\H)}
\def\fun{F_\la}
\def\sttm{\mbox{s}\tau\mbox{-tilt}\la}
\def\wte{\widetilde}
\def \text{\mbox}
\def\Mod{\mathsf{Mod}}
\hyphenation{ap-pro-xi-ma-tion}
\newcommand{\pd}{\mathsf{pd}\hspace{.01in}}
\newcommand{\add}{\mathsf{add}\hspace{.01in}}
\newcommand{\Fac}{\mathsf{Fac}\hspace{.01in}}
\newcommand{\thick}{\mathsf{thick}\hspace{.01in}}
\newcommand{\End}{\operatorname{End}\nolimits}
\renewcommand{\Mod}{\mathsf{Mod}\hspace{.01in}}
\renewcommand{\mod}{\mathsf{mod}\hspace{.01in}}
\newcommand{\proj}{\mathsf{proj}\hspace{.01in}}
\newcommand{\co}{\mathsf{Coker}\hspace{.01in}}
\newcommand{\pr}{\mathsf{pr}\hspace{.01in}}
\begin{abstract}
Let $\C$ be a triangulated category with shift functor $[1]$ and $\R$ a rigid subcategory of $\C$.
We introduce the notions of  two-term $\R[1]$-rigid subcategories, two-term (weak) $\R[1]$-cluster tilting subcategories and two-term maximal $\R[1]$-rigid subcategories, and discuss relationship between them. Our main result shows that
 there exists a bijection between the set of two-term $\R[1]$-rigid subcategories of $\C$ and the set of $\tau$-rigid subcategories of $\mod\R$, which induces a one-to-one correspondence between the set of two-term weak $\R[1]$-cluster tilting subcategories of $\C$ and the set of support $\tau$-tilting subcategories of $\mod\R$. This generalizes the main results in \cite{YZZ} where $\R$ is a cluster tilting subcategory. When $\R$ is a silting subcategory, we prove that the two-term weak $\R[1]$-cluster tilting subcategories are precisely two-term silting subcategories in \cite{IJY}. Thus the bijection above induces the bijection given by Iyama-J{\o}rgensen-Yang in \cite{IJY}.\\[0.2cm]
\textbf{Key words:} $\R[1]$-rigid subcategories; $\R[1]$-cluster tilting subcategories; $\tau$-rigid subcategories;
support $\tau$-tilting subcategories; silting subcategories; $d$-rigid subcategories. \\[0.2cm]
\textbf{ 2010 Mathematics Subject Classification:} 18E30.\end{abstract}

\medskip

\pagestyle{myheadings}
\markboth{\rightline {\scriptsize   P. Zhou and B. Zhu}}
         {\leftline{\scriptsize Two-term relative cluster tilting subcategories, $\tau-$tilting modules and silting subcategories}}
\thispagestyle{plain}

\section{Introduction}

In \cite{AIR}, Adachi, Iyama and Reiten introduced a generalization of classical tilting theory, which are called $\tau$-tilting theory.
They proved that for a $2$-Calabi-Yau triangulated category $\C$ with a basic cluster tilting object $T$, there exists a
 bijection between cluster tilting objects in $\C$ and support $\tau$-tilting modules in $\mod\End_\C(T)^{\textrm{op}}$.
Note that each cluster-tilting object is maximal rigid in a $2$-Calabi-Yau triangulated category. But
the converse is not true in general. Chang-Zhang-Zhu \cite{CZZ} and Liu-Xie \cite{LX}
showed that the similar version of the above bijection is also valid for a $2$-Calabi-Yau triangulated category with a basic rigid object. Unfortunately, many examples (see for example \cite[Example 2.15]{YZ}) indicate that the Adachi-Iyama-Reiten's bijection does not hold if  $\C$ is not $2$-Calabi-Yau.
It is then reasonable to find a class of objects in $\C$ which correspond to support $\tau$-tilting modules in $\mod\End_\C(T)^{\textrm{op}}$ bijectively in more general setting. For these purposes, Yang and Zhu \cite{YZ} introduced the notion of relative cluster tilting objects in a triangulated category $\C$, which are a generalization of cluster-tilting objects. Let $\C$ be a triangulated category and $T$ a cluster tilting object in $\C$.
They established a one-to-one correspondence between the $T[1]-$cluster tilting objects of $\C$ and the support $\tau$-tilting modules over $\End_\C(T)^{\textrm{op}}$.
This bijection is generalized  by Fu, Geng and Liu \cite{FGL} recently to rigid object. Let $\C$ be a triangulated category with shift functor $[1]$ and $R\in \C$ a basic rigid object with endomorphism algebra $\Gamma$. They introduced the notion of the $R[1]$-rigid objects in the finitely presented subcategory $\pr R$ of $\C$ and show that there exists a bijection between the set of basic $R[1]$-rigid objects in $\pr R$ and the set of basic $\tau$-rigid pairs of $\Gamma$-modules, which induces a one-to-one correspondence between the set of basic maximal $R[1]$-rigid objects with respect to $\pr R$ and the set of basic support $\tau$-tilting $\Gamma$-modules.

On the other hand, in \cite{IJY}, Iyama, J{\o}rgensen and Yang gave a functor version of $\tau$-tilting theory. They consider modules over a category and showed for a triangulated category $\C$ with a silting subcategory $\S$, there exists a bijection
the set of two-term silting subcategories of $\C$ (i.e. silting subcategories in $\S*\S[1]$) and the set of support $\tau$-tilting subcategories of $\mod\S$.
Let $\C$ be a triangulated category with a cluster tilting subcategory $\T$. Yang, Zhou and Zhu \cite{YZZ} introduced the notion of $\T[1]$-cluster tilting subcategories of $\C$. They showed that there exists a bijection between the set of $\T[1]$-cluster tilting subcategories of $\C$ and the set of support $\tau$-tilting subcategories of $\mod\T$ which generalizes the bijection in \cite{YZ} and is an analog to that in \cite{IJY}. The aim of the paper is to unify the bijections given in \cite{IJY} and in \cite{YZZ} by consider the two-term relative rigid subcategories. We also explain this generalization is much nature and has several other applications.

Let $\C$ be a triangulated category and $\R$ a rigid subcategory of $\C$. Motivated by the bijections given by Iyama-J{\o}rgensen-Yang \cite{IJY}, Yang-Zhou-Zhu \cite{YZZ} and Fu-Geng-Liu \cite{FGL}. In Section 2,  we introduce the notions of $\R[1]$-rigid subcategories, (weak) $\R[1]$-cluster tilting subcategories and  maximal $\R[1]$-rigid subcategories, and discuss connections between them, see Theorem \ref{main1}. In Section 3, we give the Bongartz completion for $\R[1]$-rigid subcategories, see Theorem \ref{h1}.

For two subcategories $\X$ and $\Y$ of a triangulated category $\C$, we denote by $\X\ast \Y$ the collection of objects in $\C$ consisting of all such $M\in \C$ with triangles $$X\longrightarrow M \longrightarrow Y \longrightarrow X[1],$$
where $X\in \X$ and $Y\in \Y$. In Section 4, We prove the following main result.

\begin{thm}\label{main} (see Theorem \ref{b3} and Theorem \ref{b4} for more details)
Let $\C$ be a triangulated category and $\R$ a rigid subcategory of $\C$.
The functor $\mathbb{H}\colon \C\to \Mod\, \R$ induces a bijection $$\Phi\colon \X\longmapsto \big(\mathbb{H}(\X),\R\cap\X[-1]\big)$$
from the first of the following sets to the second:
\begin{itemize}
\item[\emph{(I)}] two-term $\R[1]$-rigid subcategories of $\C$ ("two-term" means it is contained in $\RR$).

\item[\emph{(II)}] $\tau$-rigid pairs of $\mod\, \R$.

\end{itemize}
It restricts to a bijection from the first to the second of the following sets.
\begin{itemize}
\item[\emph{(I)}] Two-term weak $\R[1]$-cluster tilting subcategories of $\C$.

\item[\emph{(II)}] Support $\tau$-tilting subcategories of $\mod\, \R$.
\end{itemize}
\end{thm}

When $\R$ is a cluster tilting subcategory of $\C$, we show that the bijection reduces to the bijection between the set of $\R[1]$-cluster tilting subcategories of $\C$ and the set of  support $\tau$-tilting subcategories obtained by Yang-Zhou-Zhu \cite{YZZ}.

When $\R$ is a silting subcategory of $\C$, we show that the two-term weak $\R[1]$-cluster tilting subcategories of $\C$ coincide with the two-term silting subcategories of $\C$, see Section 5 and see Theorem \ref{main2}. As an application, Theorem \ref{main} recovers the bijection between the set of two-term silting subcategories of $\C$ and the set of support $\tau$-tilting subcategories of $\mod\R$ obtained by Iyama-J{\o}rgensen-Yang \cite{IJY}.

We conclude this section with some conventions. Throughout this article, $k$ is an algebraically closed field.
When we say that $\C$ is a category, we always assume that $\C$ is a Hom-finite Krull-Schmidt $k$-linear category.
Let $\C$ be an additive category. When we say that $\X$ is a subcategory of $\C$, we always assume that $\X$ is a full subcategory which is closed under isomorphisms, direct sums and direct summands.
$^\perp{\X}$ denotes the subcategory consisting of $Y\in\C$ with $\Hom_{\C}(Y,X)=0$ for any $X\in\X$, and $\X^{\perp}$
denotes the subcategory consisting of $Y\in\C$ with $\Hom_{\C}(X,Y)=0$ for any $X\in\X$.
We denote by $[\X]$ the ideal of $\C$ consisting of morphisms which factor through objects in $\X$. For any object $M$, we denote by $\add M$ the full subcategory of $\C$ consisting of direct summands of direct sum of finitely many copies of $M$. Let $\X$ and $\Y$ be subcategories of $\C$. We denote by $\X \vee \Y$ the smallest subcategory of $\C$ containing $\X$ and $\Y$.

Let $\C$ be a triangulated category with a shift functor $[1]$. For objects $X$ and $Y$ in $\C$, we define
$\text{Ext}_{\C}^i(X, Y)= \mbox{Hom}_{\C}(X, Y[i]).$
For two subcategories $\X$ and $\Y$ of $\C$, we denote by $\Ext^1_{\C}(\X,\Y)=0$ if $\Ext^1_{\C}(X,Y)=0$ for any $X\in \X, Y\in \Y$.

 Recall that a triangulated cateory $\C$ is called \emph{$2$-Calabi-Yau} if there exists a bifunctorial isomorphism \[\Hom_{\C}(X,\ Y)\simeq\D\Hom_{\C}(Y,\ X[2]) ~\text{for any}~ X,Y\in \C,\]
where $\D=\Hom_k(-,k)$ is the usual duality over $k$.

\section{Relative rigid subcategories and related subcategories}
\setcounter{equation}{0}

 An important class of subcategories of a triangulated category are the cluster tilting subcategories, which have many nice properties.  We recall the definition of cluster tilting subcategories and related subcategories from \cite{BMRRT,KR,KZ,IY}.

\begin{defn}Let $\C$ be a triangulated category.
\begin{enumerate}
\item[\emph{(i)}] A subcategory $\R$ of $\C$ is called {\rm rigid} if ${\rm Ext}_{\C}^1(\R, \R)=0$.
\item[\emph{(ii)}] A subcategory $\R$ of $\C$ is called {\rm maximal rigid} if it is rigid and maximal with respect to the property: $\R=\{M\in\C\ |\ {\rm Ext}_{\C}^1(\R\vee \add M, \R\vee \add M)=0\}$.
\item[\emph{(iii)}] A functorially finite subcategory $\R$ of $\C$ is called {\rm cluster tilting} if
$$\R=\{M\in\C\ |\ {\rm Ext}_{\C}^1(\R, M)=0\}=\{M\in\C\ |\ {\rm Ext}_{\C}^1(M, \R)=0\}.$$
\item[\emph{(iv)}] An object $R$ in $\C$ is called rigid, maximal rigid, or cluster tilting if $\add R$ is rigid, maximal rigid, or cluster tilting respectively.

\end{enumerate}
\end{defn}
\begin{rem}\label{rem:ctsubcat}
In fact, Koenig and Zhu \emph{\cite{KZ}} indicate that a subcategory $\R$ of $\C$ is cluster tilting if and only if it is contravariantly finite in $\C$ and\ \  $\R=\{M\in\C\ |\ {\rm Ext}_{\C}^1(\R, M)=0\}$.
\end{rem}

Let $\C$ be a triangulated category and $\R$ a rigid subcategory of $\C$. By \cite[Proposition 2.1(1)]{IY}, $\RR$ is closed
under direct summands. In other words, $\RR$ is a subcategory of $\C$.

In the following we introduce the notion of relative cluster tilting subcategories and related objects, compare with \cite{YZ,YZZ,FGL}.
\begin{defn}\label{a1}
Let $\C$ be a triangulated category and $\R$ a rigid subcategory of $\C$.
\begin{itemize}
\item[\emph{(i)}] A subcategory $\X$ in $\C$ is called  \emph{$\R[1]-$rigid }if  $[\R[1]](\X, \X[1])=0$. Any $\R[1]-$rigid subcategory in $\RR$ is called two-term $\R[1]-$rigid.

\item[\emph{(ii)}] A subcategory $\X\subseteq \RR$ is called \emph{two-term $\R[1]-$maximal rigid} if $\X$ is $\R[1]$-rigid and for any $M\in \RR$,
$$[\R[1]](\X\vee\add M, (\X\vee\add M)[1])=0 \text{ implies } M\in\X. $$

\item[\emph{(iii)}] A subcategory $\X\subseteq \RR$ is called {\rm two-term weak $\R[1]-$cluster tilting} if $\R\subseteq\X[-1]\ast\X$ and  $$\X=\{M\in\RR\ |\ [\R[1]](M, \X[1])=0 \;\, \emph{and } \ [\R[1]](\X, M[1])=0 \ \}.$$

\item[\emph{(iv)}] A subcategory $\X\subseteq \RR$ is called {\rm two-term $\R[1]-$ cluster tilting} if $\X$ is contravariantly finite and  $$\X=\{M\in\RR\ |\ [\R[1]](M, \X[1])=0 \;\, \emph{and } \ [\R[1]](\X, M[1])=0 \ \}.$$

\item[\emph{(v)}] An object $X$ is called two-term $\R[1]$-rigid, two-term maximal $\R[1]$-rigid, two-term weak $\R[1]$-cluster tilting, or two-term $\R[1]$-cluster tilting if $\add X$ is two-term $\R[1]$-rigid, two-term maximal $\R[1]$-rigid, two-term weak $\R[1]$-cluster tilting, or two-term $\R[1]$-cluster tilting respectively.
\end{itemize}
\end{defn}

When $\R$ is a cluster tilting subcategory, then $\RR=\C$, the notations above are usual ones studied in \cite{YZZ},\cite{YZ}. It is well-known that cluster tilting subcategories are functorially finite maximal rigid, but the converse is not true in general. From definition, one can easily see two-term $\R[1]-$cluster tilting are contravariantly finite two-term maximal $\R[1]-$rigid.
The main result of this section is the converse also holds, i.e. any contravariantly finite two-term maximal $\R[1]$-rigid subcategory is two-term $\R[1]$-cluster tilting. This is a generalization of Theorem 2.6 in \cite{ZZ}.
From now on to the end of this article, we assume that $\C$ is a triangulated category with a rigid subcategory $\R$.

\begin{thm}\label{main1}
Any two-term $\R[1]$-cluster tilting subcategory is  precisely contravariantly finite two-term maximal $\R[1]$-rigid subcategory.
\end{thm}

As a direct consequence, we have the following important result.

\begin{cor}\label{a4}{\emph{\cite[Theorem 2.6]{ZZ}}}
Let $\C$ be a $2$-Calabi-Yau triangulated category with a cluster tilting subcategory $\cal R$. Then every
functorially finite maximal rigid subcategory is cluster-tilting.
\end{cor}
\proof For $2-$Calabi-Yau triangulated category $\C$, $\X$ is $\R[1]-$rigid if and only if it is rigid, and since $\R*\R[1]=\C$, any functorially finite maximal rigid subcategory is functorially finite two-term maximal $\R[1]-$rigid. Thus it is cluster tilting subcategory following from Theorem \ref{main1}.  \qed

In order to prove Theorem \ref{main1}, we need the following lemma.

\begin{lem}\label{keylem}
Let $\X$ be a two-term maximal $\R[1]$-rigid subcategory.
\begin{itemize}
\item[\emph{(a)}]For any object $R_0\in\R$, if there exists a triangle: $$M[-1]\s{f}\longrightarrow R_0 \s{g}\longrightarrow X \s{h}\longrightarrow M$$ such that $g\colon R_0\longrightarrow X$ is a left $\X$-approximation of $R_0$, then $M\in\X$.
\item[\emph{(b)}] For any object $R_0\in\R$, if there exists a triangle: $$M[-1]\s{f}\longrightarrow X[-1] \s{g}\longrightarrow R_0 \s{h}\longrightarrow M$$ such that $g\colon X[-1] \longrightarrow R_0$ is a right $\X[-1]$-approximation of $R_0$, then $M\in\X$.
\end{itemize}
\end{lem}

\proof We only prove (a), the proof of (b) is similar. For $R_0\in \R$, suppose we have the triangle $$M[-1]\s{f}\longrightarrow R_0 \s{g}\longrightarrow X \s{h}\longrightarrow M$$ such that $g\colon R_0\longrightarrow X$ is a left $\X$-approximation of $R_0$. Object $X\in\X\subseteq\RR$, there exists a triangle
$$R_1\s{~}\longrightarrow R_2 \s{~}\longrightarrow X \s{~}\longrightarrow R_1[1]$$
where $R_1,R_2\in\R$. By the octahedral axiom,
we have the following commutative diagram
\begin{center}
$\xymatrix@!@C=0.5cm@R=0.5cm{
               &R_2\ar[d]\ar@{=}[r]  & R_2\ar[d]  \\
R_0\ar@{=}[d]\ar[r]  & X\ar[d]\ar[r]     & M\ar[r]\ar[d] & R_0[1]\ar@{=}[d]  \\
R_0\ar[r]^{h_1\;\;}           & R_1[1]\ar[d]\ar[r]    & N\ar[d]\ar[r]    & R_0[1] \\
               & R_2[1]\ar@{=}[r]      & R_2[1]
  }$
\end{center}
of triangles. Since $\R$ is a rigid subcategory, we have $h_1=0$. It follow that $N\simeq R_1[1]\oplus R_0[1]\in\R[1]$
implies $M\in\RR$.

For any $ x\in [\R](M[-1], X_1),$ where $X_1\in \X$, there are two morphisms $x_1\colon M[-1]\rightarrow R_1$ and $x_2\colon R_1\rightarrow X_1$ such that $x=x_2x_1$, where $R_1\in \R$.

\begin{center}
$\xymatrix@!@C=0.8cm@R=0.4cm{
X_0[-1]\ar[r]^{h[-1]} & M[-1]\ar[r]^{f}\ar[d]^{x_1} & R_0\ar[r]^{g}\ar@{.>}[ddl]^{a} & X\ar[r]^{h}\ar@{.>}[ddll]^{b} & M\ar[r]^{-f[1]\quad} & R_0[1] \\
                       & R_1\ar[d]^{x_2}            &               &               & R_2\ar[u]^{y_2}\ar@{.>}[lu]^{c}    \\
                      & X_1                         &                &               & X_2[-1]\ar[u]^{y_1}
  }$

\medskip
\end{center}
Since $\X$ is $\R[1]$-rigid, we have $xh[-1]=x_2(x_1h[-1])=0$, then there exists $a\colon R_0\to X_1$ such that $x=af$. Since $g$ is a left $\X$-approximation of $R_0$, we know that there exists a morphism $b\colon X\to X_1$ such that $a=bg$. Therefore, $x=af=b(gf)=0$ and
\begin{equation}\label{equi:lem1}
\begin{array}{l}
[\R[1]](M, \X[1])=0.
\end{array}
\end{equation}
  For any  $ y\in [\R](X_2[-1], M),$ where $X_2\in \X$, there are two morphisms $y_1\colon X_2[-1]\rightarrow R_2$ and $y_2\colon R_2\rightarrow M$ such that $y=y_2y_1$, where $R_2\in \R$ and $X_2\in\X$. Since $f[1]y_2=0$, there exists $c\colon R_2\to X$ such that $y_2=hc$. Since $\X$ is $\R[1]$-rigid, we have $y=y_2y_1=h(cy_1)=0.$ Therefore,
\begin{equation}\label{equi:lem2}
\begin{array}{l}
[\R[1]](\X, M[1])=0.
\end{array}
\end{equation}
For any $z\in [\R](M[-1], M),$ there are two morphisms $z_1\colon M[-1]\rightarrow R_3$ and $z_2\colon R_3\rightarrow M$ such that $z=z_2z_1$, where $R_3\in \R$. Since $f[1]z_2=0$, there exists $d\colon R_3\to X$ such that $z_2=hd$.
\begin{center}
$\xymatrix@!@C=0.5cm@R=0.3cm{
 R_0\ar[r]^{g} & X\ar[r]^{h} & M\ar[r]^{-f[1]\quad} & R_0[1] \\
 &               & R_3\ar[u]^{z_2}\ar@{.>}[lu]^d    \\
&               & M[-1]\ar[u]^{z_1}
  }$
\end{center}
By equality (\ref{equi:lem1}), we have $z=z_2z_1=h(dz_1)=0$. Thus,
\begin{equation}\label{equi:lem3}
\begin{array}{l}
[\R[1]](M, M[1])=0.
\end{array}
\end{equation}
Using equalities (\ref{equi:lem1}), (\ref{equi:lem2}) and (\ref{equi:lem3}), we get $[\R[1]](\X\vee\add M, (\X\vee\add M)[1])=0$. Since $\X$ is maximal $\R[1]$-rigid, we have $M\in \X$.  This concludes the proof of (a). \qed
\medskip

The following corollary is a direct consequence of Lemma \ref{keylem}.

\begin{cor}\label{cor1}
Let $\X$ be a covariantly finite or contravariantly finite two-term maximal $\R[1]$-rigid subcategory. Then $\R\subseteq\X[-1]\ast\X$.
\end{cor}

\textbf{Now we give the proof of Theorem \ref{main1}.}
\medskip

 It is easy to see that any two-term $\R[1]$-cluster tilting subcategory is contravariantlty finite two-term maximal $\R[1]$-rigid subcategory. We prove the other direction.

  Assume that $\X$ is a contravariantly finite two-term maximal $\R[1]$-rigid subcategory.
Clearly, $$\X\subseteq \{M\in\RR\ |\ [\R[1]](\X, M[1])=0= [\R[1]](M, \X[1])\}.$$
For any object $M\in \{M\in\RR\ |\ [\R[1]](\X, M[1])=0= [\R[1]](M, \X[1])\}$,
there exists a triangle
$$\xymatrix{R_1\ar[r]^{f}&R_0\ar[r]^{g}&M\ar[r]^{h\;\;}&R_1[1]},$$
where $R_0,R_1\in\R$. By Corollary \ref{cor1}, there exists a triangle
$$\xymatrix{R_0\ar[r]^{u}&X_1\ar[r]^{v}&X_2\ar[r]^{w\;\;}&R_0[1]},$$
where $X_1,X_2\in\X$. Since $\X$ is $\R[1]$-rigid, we have that $u$ is a left $\X$-approximation of $R_0$.
By the octahedral axiom, we have a commutative diagram
$$\xymatrix{
R_1\ar[r]^{f}\ar@{=}[d]&R_0\ar[r]^{g}\ar[d]^{u}&M\ar[r]^{h}\ar[d]^{a}&R_1[1]\ar@{=}[d]\\
R_1\ar[r]^{x=uf}&X_1\ar[r]^{y}\ar[d]^{v}&N\ar[r]^{z}\ar[d]^{b}&R_1[1]\\
&X_2\ar@{=}[r]\ar[d]^{w}&X_2\ar[d]^{c}\\
&R_0[1]\ar[r]^{g[1]}&M[1]}$$
of triangles. We claim that $x$ is a left $\X$-approximation of $R_1$.
Indeed, for any
morphism $\alpha\colon R_1\to X'$, where $X'\in\X$, since $\alpha\circ h[-1]\in[\R](M[-1],X')=0$, there exists a morphism $\beta\colon R_0\to X'$ such that $\alpha=\beta f$. Since $u$ is a left $\X$-approximation of $R_0$ and $X'\in\X$,  there exists a morphism $\gamma\colon X_1\to X'$ such that
$\beta=\gamma u$ and then $\alpha=\gamma(uf)=\gamma x$. This shows that $x$ is a left $\X$-approximation of $R_1$.

By Lemma \ref{keylem}, we have $N\in\X$. Since $$c=g[1]w\in [\R[1]](X_2,M[1])=0.$$
This shows that the triangle
$$\xymatrix{M\ar[r]^{a}&N\ar[r]^{b}&X_2\ar[r]^{c\;\;}&M[1].}$$
splits. It follows that $M$ is a direct summand of $N$ and then $M\in\X$.

Therefore, $\X$ is $\R[1]$-cluster tilting.  \qed
\medskip

\medskip

\section{Bongartz completion for relative rigid subcategories}
\setcounter{equation}{0}

The main result of this section is the following analog of Bongartz completion for relative rigid subcategories.

\begin{thm}\label{h1}
Let $\C$ be a triangulated category and $\R$ a rigid subcategory of $\C$.
If $\X$ is a contravariantly finite two-term $\R[1]$-rigid subcategory of $\C$, then
$\C$ has a two-term weak $\R[1]$-cluster tilting subcategory $\Y$ which contains $\X$.
\end{thm}

\proof For any object $R\in\R$, we take a triangle
$$R\s{f}\longrightarrow U_R \s{g}\longrightarrow X_1 \s{h}\longrightarrow R[1],$$
where $h$ is a  right $\X$-approximation of $R[1]$. Let $\Y\coloneqq \X\vee\add\{\, U_R
\,|\, R \in \R \,\}$ be the additive closure of $\X$ and $\{\, U_R
\,|\, R \in \R \,\}$. We claim that $\Y$ is a two-term weak $\R[1]$-cluster tilting subcategory of $\C$ which contains $\X$.
\smallskip

Since $X_1\in\X\subseteq\RR$, there exists a triangle
$$R''\s{~}\longrightarrow R' \s{~}\longrightarrow X_1 \s{~}\longrightarrow R''[1]$$
where $R',R''\in\R$. By the octahedral axiom,
we have the following commutative diagram
\begin{center}
$\xymatrix@!@C=0.5cm@R=0.5cm{
               &R''\ar[d]\ar@{=}[r]  & R''\ar[d]  \\
R\ar@{=}[d]\ar[r]  & Q\ar[d]\ar[r]     & R'\ar[r]^{q=0}\ar[d] & R[1]\ar@{=}[d]  \\
R\ar[r]         & U_R\ar[d]\ar[r]    & X_1\ar[d]\ar[r]    & R[1] \\
               & R''[1]\ar@{=}[r]      & R''[1]
  }$
\end{center}
of triangles. Since $\R$ is a rigid, we have $q=0$. It follow that $Q\simeq R\oplus R'\in\R$
implies $U_R\in\RR$.  Hence we obtain $\Y\subseteq\RR$.

It is clear that $\R\subseteq\Y[-1]\ast\Y$. \   It remains to show that
$$\Y=\{M\in\RR\ |\ [\R[1]](M, \Y[1])=0=[\R[1]](\Y, M[1])\}.$$

\textbf{Step 1:} We show that $\Y$ is a $\R[1]$-rigid subcategory.
\smallskip

Take any morphism $a\in [\R[1]](X, U_R[1])$, where $X\in\X$. Since $a$ factors through an object in $\R[1]$ and $\X$ is $\R[1]$-rigid, we have $g[1]\circ a=0$.
$$\xymatrix@!@C=0.8cm@R=0.6cm{
 U_R\ar[r]^g & X_1\ar[r]^h    & R[1]\ar[r]^{-f[1]\;} & U_R[1]\ar[r]^{-g[1]\;}  &  X_1[1]   \\
   &&&X\ar[u]_{a}\ar@{.>}[ul]_{b}\ar@{.>}[ull]^{c}              }
$$
Thus there exists a morphism $b\colon X\rightarrow R[1]$ such that $a=-f[1]\circ b$.
Since $h$ is a right $\X$-approximation of $R[1]$,  there exists a morphism $c\colon X\rightarrow X_1$ such that $b=hc$. It follows that $$a=-f[1]b=(-f[1]h)c=0,$$ and therefore
\begin{equation}\label{m1}
\begin{array}{l}
[\R[1]](X, U_R[1])=0.
\end{array}
\end{equation}

For any morphism $u\in [\R[1]](U_R, X[1])$, where $X\in\X$, we know that there are two morphisms $u_1\colon U_R\rightarrow R_0[1]$ and $u_2\colon R_0[1]\rightarrow X[1]$ such that $u=u_2u_1$, where $R_0\in\R$.
$$\xymatrix@!@C=0.6cm@R=0.6cm{
X_1[-1]\ar[r]^{\quad h[-1]}&R\ar[r]^f&U_R\ar[d]^{u_1} \ar[r]^g & X_1\ar@{.>}[ld]^{v}   \\
   && R_0[1]\ar[d]^{u_2}\\
   && X[1]             }
$$
Since $\Hom_{\C}(R,R_0[1])=0$, there exists a morphism $v\colon X_1\to T_0[1]$ such that $u_1=vg$.
Since $\X$ is $\R[1]$-rigid, we have $u=u_2u_1=(u_2v)g=0$. Therefore
\begin{equation}\label{m2}
\begin{array}{l}
[\R[1]](U_R, X[1])=0.
\end{array}
\end{equation}

For any morphism $x\in [\R[1]](U_R, U_R[1])$, we know that there are two morphisms $x_1\colon U_R\rightarrow R_1[1]$ and $x_2\colon R_1[1]\rightarrow U_R[1]$ such that $x=x_2x_1$, where $R_1\in \R$.
$$\xymatrix@!@C=0.6cm@R=0.6cm{
R\ar[r]^{f} & U_R\ar[r]^g\ar[d]^{x_1} & X_1\ar[r]^h\ar@{-->}[ld]^{y}& R[1]\\
   & R_1[1]\ar[d]^{x_2}\\
   & U_R[1]             }$$
Since $\Hom_{\C}(R,R_1[1])=0$, there exists a morphism $y\colon X_1\to R_1[1]$ such that $x_1=yg$.
Since $x_2y\in[\R[1]](X_1,U_R[1])$ and $[\R[1]](\X,U_R[1])=0$, we have $x=x_2x_1=(x_2y)g=0$. Therefore
\begin{equation}\label{m3}
\begin{array}{l}
[\R[1]](U_R, U_R[1])=0.
\end{array}
\end{equation}
Using equalities (\ref{m1}), (\ref{m2}) and (\ref{m3}),
we get that $\Y$ is a $\R[1]$-rigid subcategory.
\medskip

\textbf{Step 2:} We show that $\{M\in\RR\ |\ [\R[1]](M, \Y[1])=0=[\R[1]](\Y, M[1])\}\subseteq\Y.$
\medskip

For any object $M\in\RR$, assume that $[\R[1]](M, \Y[1])=0=[\R[1]](\Y, M[1])$. Then there exists a triangle
$$R_3\s{u}\longrightarrow R_2 \s{v}\longrightarrow M \s{w}\longrightarrow R_3[1],$$
where $R_3,R_2\in\R$.  For the object $R_3\in\R$, we take a triangle
$$R_3\s{\alpha}\longrightarrow U_{R_3} \s{\beta}\longrightarrow X_3 \s{\gamma}\longrightarrow R_3[1],$$
where $\gamma$ is a  right $\X$-approximation of $R_3[1]$ and $U_{R_3}\in\Y$.
\medskip

We claim that $-u[1]\circ\gamma\colon X_3\to R_2[1]$ is a right $\X$-approximation of $R_2[1]$. Indeed, for any morphism $s\colon X\to R_2[1]$ with $X\in\X$. Since $v[1]\circ s\in [\R[1]](\X,M[1])\subseteq[\R[1]](\Y,M[1])=0$, there exists a morphism $t\colon X\to X_3$ such that $s=-u[1]\circ t$.
$$\xymatrix@!@C=0.8cm@R=0.5cm{
&&&&X\ar[d]^{s}\ar@{.>}[dl]_{t}\\
R_2\ar[r]^{u}& R_3\ar[r]^{v} & M\ar[r]^{w}& R_3[1]\ar[r]^{-u[1]}&R_2[1]\ar[r]^{-v[1]}&M[1].}$$
Since $\gamma$ is a right $\X$-approximation of $R_3[1]$ and $X\in\X$, there exists a morphism
$\phi\colon X\to X_3$ such that $t=\gamma\phi$. It follows that $s=-u[1]\circ t=(-u[1]\circ\gamma)\phi$.
This show that $-u[1]\circ\gamma\colon X_3\to R_2[1]$ is a right $\X$-approximation of $R_2[1]$. Thus there exists a
triangle
$$R_2\xrightarrow{~d~}U_{R_2}\xrightarrow{~e~}X_3\xrightarrow{-u[1]\circ\gamma}R_2[1]$$
where $-u[1]\circ\gamma$ is a right $\X$-approximation of $R_2[1]$ and $U_{R_2}\in\Y$.
\medskip

By the octahedral axiom, we have a commutative diagram
\begin{center}
$\xymatrix@!@C=0.8cm@R=0.5cm{
               &M\ar[d]^{w}\ar@{=}[r]  & M\ar[d]^{a} \\
X_3\ar@{=}[d]\ar[r]^{\gamma}  & R_3[1]\ar[d]^{-u[1]}\ar[r]^{-\alpha[1]\;}    & U_{R_3}[1]\ar[r]^{-\beta[1]\;\,}\ar[d]^{b} & X_3[1]\ar@{=}[d]  \\
X_3\ar[r]^{-u[1]\circ\gamma\;\;\;}           & R_2[1]\ar[d]^{-v[1]}\ar[r]^{\delta}     & N[1]\ar[d]^c\ar[r]^{\psi\;}    & X_3[1] \\
               & M[1]\ar@{=}[r]      & M[1]
  }$
\end{center}
of triangles. Since $a=-\alpha[1]\circ w\in[\R[1]](M,\R[1])=0$, we have that
$$M\xrightarrow{~a~}U_{R_3}[1]\xrightarrow{~b~}N[1]\xrightarrow{~c~}M[1]$$
splits. This implies that $M$ is a direct summand of $N$.
\medskip

Consider the following
commutative diagram
$$\xymatrix@C=1.2cm{R_2\ar[r]^{d}\ar@{=}[d]&U_{R_2}\ar[r]^{e}\ar@{-->}[d]^{\eta}&X_3\ar[r]^{-u[1]\circ\gamma\quad}\ar@{=}[d]&R_2[1]\ar@{=}[d]\\
R_3\ar[r]^{\delta[-1]}&N\ar[r]^{\psi[-1]}&X_3\ar[r]^{-u[1]\circ\gamma\quad}&R_2[1],}$$
there exists an isomorphism $\eta\colon U_{T_2}\to N$ which makes
the above diagram commutative. Thus we have $N\simeq U_{R_2}\in\Y$ and then $M\in\Y$.
This shows that $$\{M\in\RR\ |\ [\R[1]](M, \Y[1])=0=[\R[1]](\Y, M[1])\}\subseteq\Y.$$
Hence $\Y=\X\vee\add\{\, U_R\,|\, R \in \R \,\}$  is a weak $\R[1]$-cluster tilting subcategory of $\C$. \qed

\begin{lem}\label{h3}
Let $\C$ be an additive category and $\X\subseteq\cal A$ two subcategories of $\C$.
If $\X$ is contravariantly finite in $\C$, then $\cal A$ is contravariantly finite in $\C$
if and only if $\cal A/\X$ is contravariantly finite in $\C/\X$.
\end{lem}

\proof  Assume that $\cal A/\X$ is contravariantly finite in $\C/\X$. For any $C\in\C$, let $\overline{f}\colon A\to C$
be a right $\cal A/\X$-approximation of $C$. Since $\X$ is contravariantly finite in $\C$, we can take a right $\X$-approximation $\alpha\colon X\to C$. We claim that $(f,\alpha)\colon A\oplus X\to C$ is a right $\cal A$-approximation of $C$.
Indeed, for any morphism $a\colon A'\to C$, where $A'\in\cal A$, since $\overline{f}$ is a right $\cal A/\X$-approximation, there
exists a morphism $\overline{b}\colon A'\to A$ such that $\overline{f}\circ\overline {b}=\overline {a}$.
Then $a-fb$ factors through $\X$, i.e., there exist morphisms $s\colon A'\to X_0$ and $t\colon X_0\to C$, where $X_0\in\X$ such
that $a-fb=ts$. Since $\alpha$ is a right $\X$-approximation and $X_0\in\X$, there exists a morphism $c\colon X_0\to X$
such that $t=\alpha c$. It follows that $a=fb+\alpha cs=(f,\alpha)\binom{b}{cs}$. This shows that $(f,\alpha)\colon A\oplus X\to C$ is a right $\cal A$-approximation of $C$. Therefore, $\cal A$ is contravariantly finite in $\C$.
\medskip

Conversely, if $f\colon A\to C$
is a right $\cal A$-approximation of $C$, then it is easy to see that $\overline{f}\colon A\to C$
is a right $\cal A/\X$-approximation of $C$. It follows that if $\cal A$ is contravariantly finite in $\C$,
then $\cal A/\X$ is contravariantly finite in $\C/\X$.  \qed
\medskip

The above theorem immediately yields the following important conclusion.

\begin{cor}\label{h4}{\emph{\cite[Theorem 1.6]{AV}}}
Let $\C$ be a $2$-Calabi-Yau triangulated category with a cluster-tilting subcategory $\R$, and $\X$ a functorially finite rigid subcategory of $\C$. Then
$\C$ has a cluster-tilting subcategory $\Y$ which contains $\X$.
\end{cor}

\proof By Theorem \ref{h1}, we can obtain that $\Y$ exists through this way.
We now show that $\Y$ is a contravariantly finite subcategory of $\C$.
\smallskip

Since $\X$ is a covariantly finite rigid subcategory of $\C$, we have that $({^\perp{\X[1]}}, \X)$ is a cotorsion pair on $\C$. It follows that ${^\perp{\X[1]}}$ is a contravariantly finite subcategory of $\C$.
\smallskip

It suffices to show that $\Y$ is a contravariantly finite subcategory of ${^\perp{\X[1]}}$. For any object $M\in{^\perp{\X[1]}}=\X[-1]^{\perp}$, since $\R$ is a cluster-tilting, we take a right approximation $u\colon R_0\to M$ of $M$. For the object $R_0\in\R$, we take a triangle
$$R_0\s{f}\longrightarrow U_{R_0} \s{g}\longrightarrow X_0 \s{h}\longrightarrow R_0[1],$$
where $h$ is a  right $\X$-approximation of $R_3[1]$ and $U_{R_0}\in\Y$. Since $M\in{^\perp{\X[1]}}=\X[-1]^{\perp}$,
there exists a morphism $v\colon U_{R_0}\to M$ such that $u=vf$.
$$\xymatrix{
    X_0[-1]\ar[r]^{\;\; h[-1]}&R_0 \ar[r]^{f} \ar[d]^{u} & U_{R_0}\ar@{.>}[dl]^{v} \ar[r]^{g} &X_0\ar[r]^{h}&R_0[1]\\
    &M&}$$
We claim that $\overline{v}\colon U_{R_0}\to M$ is a right $\overline{\Y}$-approximation of $M$. Indeed,
let $\overline{a}\colon U_{R_1}\to M$ be any morphism, where $U_{R_1}\in\overline{\Y}$ and $R_1\in\R$.
\medskip

For the object $R_1\in\R$, there exists a triangle
$$R_1\s{x}\longrightarrow U_{R_1} \s{y}\longrightarrow X_1 \s{z}\longrightarrow R_1[1],$$
where $z$ is a  right $\X$-approximation of $R_1[1]$ and $U_{R_1}\in\R$.
Since $u\colon R_0\to M$ is a right $\R$-approximation of $M$, there exists a morphism $b\colon R_1\to R_0$ such that
$ub=ax$. Since $h$ is a  right $\X$-approximation of $R_3[1]$, there exists a morphism $c\colon X_1\to X_0$ such that
$hc=b[1]\circ z$. Thus we have the following commutative diagram
$$\xymatrix{R_1\ar[r]^{x}\ar@{=}[d]&U_{R_1}\ar[r]^{y}\ar@{.>}[d]^{d}&X_1\ar[r]^{z\;\,}\ar[d]^{c}&R_1[1]\ar[d]^{b[1]}\\
R_0\ar[r]^{f}&U_{R_0}\ar[r]^{g}&X_0\ar[r]^{h\;\,}&R_1[1],}$$
of triangles. It follows that $ax=ub=vfb=vdx$ and then $(a-vd)x=0$. Then there exists a morphism
$e\colon X_1\to M$ such that $a-vd=ey$ and then $\overline{a}=\overline{v}\circ\overline{d}$.
\medskip

This shows that $\overline{v}\colon U_{R_0}\to M$ is a right $\overline{\Y}$-approximation of $M$.
\medskip

Hence $\Y/\X$ is a contravariantly finite subcategory of ${^\perp{\X[1]}}/\X$.
By Lemma \ref{h3}, we have that $\Y$ is a contravariantly finite subcategory of ${^\perp{\X[1]}}$.
\medskip

This follows from Theorem \ref{h1}. \qed

\section{Relative rigid subcategories and $\tau$-rigid subcategories }
\setcounter{equation}{0}

Let $\C$ be a triangulated category and $\R$ a rigid subcategory of $\C$. We write $\Mod\R$ for the abelian group of contravariantly additive functor from $\R$ to the category of abelian group, and $\mod\R$ for the full subcategory of
finitely presentation functor, see \cite{Au}. There exists a functor
\begin{eqnarray*}
\mathbb{H}\colon\C & \longrightarrow & \Mod\R \\
 M& \longmapsto & \Hom_{\C}(-,M)\mid_{\R}
\end{eqnarray*}
sometimes known as the restricted Yoneda functor.

\begin{thm}\label{b1}
\begin{itemize}
\item[\emph{(i)}] \emph{\cite{Au}} For
    $M \in \Mod\, \T$ and $R \in \R$, there exists a
    natural isomorphism
\[
 \emph{\Hom}_{ \Mod\, \R }\big( \R( -,R ) , M \big)
  \xrightarrow{~\sim~} M( R ).
\]

  \item[\emph{(ii)}] \emph{\cite[Proposition 6.2]{IY}} The functor \ $\mathbb{H}$
    induces an equivalence
$$(\RR) / [ \R[1] ]\xrightarrow{~\sim~}  \mod\, \R,$$
\end{itemize}
\end{thm}

\begin{defn}{\emph{\cite[Definition 1.3]{IJY}}}
Let $\R$ be an essentially small additive category.
\begin{itemize}
\item[\emph{(i)}] Let $\M$ be a subcategory of $\mod\,\R$.  A class $\{\,
  P_1 \stackrel{ \pi^M }{ \rightarrow } P_0 \rightarrow M \rightarrow
  0 \,\mid\, M \in \M \,\}$ of projective presentations in $\mod\,
  \R$ is said to have \emph{Property (S)} if
\[
  \emph{\Hom}_{ \mod\,\R }( \pi^M , M' )
  : \emph{\Hom}_{ \mod\,\R }( P_0 , M' )
    \rightarrow \emph{\Hom}_{ \mod\, \R }( P_1 , M' )
\]
is surjective for any $M , M' \in \M$.

\item[\emph{(ii)}] A subcategory $\M$ of $\mod\, \R$ is said to be
  \emph{$\tau$-rigid} if there is a class of projective presentations
  $\{P_1 \rightarrow P_0 \rightarrow M \rightarrow 0\mid M\in\M\}$
  which has Property (S).

\item[\emph{(iii)}] A \emph{$\tau$-rigid pair} of $\mod\, \R$ is a pair $(
  \M , \E )$, where $\M$ is a $\tau$-rigid subcategory of $\mod\,
  \R$ and $\E \subseteq \R$ is a subcategory with $\M( \E )=0$,
  that is, $M( E )=0$ for each $M \in \M$ and $E \in \E$.

\item[\emph{(iv)}] A $\tau$-rigid pair $( \M , \E )$ is \emph{support
  $\tau$-tilting} if $\E = \emph{Ker}\,( \M )$ and for each $R\in \R$
  there exists an exact sequence $\R( - , R )
  \stackrel{f}{\rightarrow} M^0 \rightarrow M^1 \rightarrow 0$ with
  $M^0, M^1 \in \M$ such that $f$ is a left $\M$-approximation. In this case, $\M$ is called a \emph{support
  $\tau$-tilting subcategory} of $\mod\R$.
\end{itemize}
\end{defn}

\begin{lem}\label{b2}
Let $\X$ be a subcategory of $\C$.
For any object $X\in\X$, let
\begin{equation}\label{m1}
\begin{array}{l}
R_1\s{f}\longrightarrow R_0 \s{g}\longrightarrow X \s{h}\longrightarrow R_1[1]
\end{array}
\end{equation}
be a triangle in $\C$ with $R_0,R_1\in\R$. Then the functor $\mathbb{H}$ gives a projective presentation
\begin{equation}\label{m2}
\begin{array}{l}
P_1^{\mathbb{H}(X)}\xrightarrow{~~~\pi^{\mathbb{H}(X)}}P_0^{\mathbb{H}(X)}\xrightarrow{~~}\mathbb{H}(X)\xrightarrow{~~}0
\end{array}
\end{equation}
in $\mod\, \R$, and $\X$ is a $\R[1]$-rigid if and only if the class $\{\, \pi^{\mathbb{H}(X)} \,|\, X \in \X \,\}$ has Property (S).
\end{lem}

\proof It is clear that $\mathbb{H}$ applies to the triangle (\ref{m1}) gives the projective
presentation (\ref{m2}).\\[0.2cm]
By Theorem \ref{b1}(i), the morphism
$\Hom_{\mod\,\T}\big(\pi^{\mathbb{H}(X)},\mathbb{H}(X')\big)$, where $X'\in\X$ is the same as
\begin{equation}\label{t3}
\begin{array}{l}
\Hom_{\C}( R_0, X' )\xrightarrow{\Hom_{\C}(f,\ X')} \Hom_{\C}(R_1,X').
\end{array}
\end{equation}
Therefore, the class $\{\, \pi^{\mathbb{H}(X)} \,|\, X\in \X \,\}$ has Property (S) if and
only if the morphism (\ref{t3}) is
surjective for all $X, X'\in\X$.
\vspace{1mm}

Assume that the class $\{\, \pi^{\mathbb{H}(X)} \,|\, X\in \X \,\}$ has Property (S). For any morphism
$a\in[\R[1]](\X,\X[1])$, we know that there exist two morphisms
$a_1\colon X\to R[1]$ and $a_2\colon R[1]\to X'[1]$ such that $a=a_2a_1$, where $X,X'\in\X$ and  $R\in\R$.
Since $\Hom_{\C}(R_0,R[1])=0$, there exists a morphism $b\colon R_1[1]\to R[1]$ such that
$a_1=bh$.
$$\xymatrix@!@C=0.5cm@R=0.3cm{
 R_1\ar[r]^{f} & R_0\ar[r]^{g} & X\ar[r]^{h\quad} \ar[d]^{a_1}& R_1[1]\ar@{.>}[dl]^{b} \\
 &               & R[1]\ar[d]^{a_2} \\
&               & X'[1] }$$
Since $\Hom_{\C}(f,\ X')$ is surjective, there exists a morphism $c\colon R_0\to X'$ such that
$a_2[-1]\circ b[-1]=cf$ and then $a_2b=c[1]\circ f[1]$. It follows that
$a=a_2a_1=a_2bh=c[1]\circ(f[1]h)=0$.

This shows that
$[\R[1]](\X,\X[1])=0$. Hence $\X$ is a $\R[1]$-rigid.
\smallskip

Conversely, assume that $\X$ is a $\R[1]$-rigid subcategory. For any morphism
$x\colon R_1\to X'$, since $\X$ is $\R[1]$-rigid, we have $x\circ h[-1]=0$. Then there exists a morphism
$y\colon R_0\to X'$ such that $x=yf$.
$$\xymatrix@!@C=0.6cm@R=0.3cm{
 X[-1]\ar[r]^{\; h[-1]}&R_1\ar[r]^{f}\ar[d]^{x}& R_0\ar@{.>}[dl]^{y}\ar[r]^{g} & X\ar[r]^{h\quad}& R_1[1]\\
 &X'&& }$$
This shows that $\Hom_{\C}(f, X')\colon \Hom_{\C}( R_0, X' )\to\Hom_{\C}(R_1,X')$ is surjective. By the above discussion, we obtain that the class $\{\, \pi^{\mathbb{H}(X)} \,|\, X\in \X \,\}$ has Property (S). \qed

\begin{thm}\label{b3}
Let $\C$ be a triangulated category and $\R$ a rigid subcategory of $\C$.
The functor $$\mathbb{H}\colon \C\to \Mod\, \R$$ induces a bijection $$\Phi\colon \X\longmapsto \big(\mathbb{H}(\X),\R\cap\X[-1]\big)$$
from the first of the following sets to the second:
\begin{itemize}
\item[\emph{(I)}] $\R[1]$-rigid subcategories of ~$\C$ which are contained in $\RR$.

\item[\emph{(II)}] $\tau$-rigid pairs of $\mod\, \R$.
\end{itemize}
\end{thm}

\proof \textbf{Step 1:} The map $\Phi$ has values in $\tau$-rigid pairs of $\mod\, \T$.
\smallskip

Assume that  $\X$ is a $\R[1]$-rigid subcategories of ~$\C$ which are contained in $\RR$.
Then for any $X\in\X$, there exists a triangle in $\C$
$$R_1\s{f}\longrightarrow R_0 \s{g}\longrightarrow X \s{h}\longrightarrow R_1[1],$$
where $R_0,R_1\in\R$.  By Lemma \ref{b2}, we have that
$\mathbb{H}$ sends the set of these triangles to a set of  projective presentations
(\ref{m2}) which has Property (S).

It remains to show that for any $X\in\X$ and $X'\in\R\cap\X[-1]$, we have $\mathbb{H}(X)(X')=0$. Indeed,
since $\X$ is a $\R[1]$-rigid, we have $\mathbb{H}(X)(X')=\Hom_{\C}(X',X)=0$.
$$\xymatrix{
X'\in\X[-1]\ar[r]\ar@{=}[dr]&X\in\X\\
 &X'\in\R\ar[u] }$$
This shows that $\big(\mathbb{H}(\X),\R\cap\X[-1]\big)$ is a $\tau$-rigid pair of $\mod\, \R$.
\medskip

\textbf{Step 2:} The map $\Phi$ is surjective.
\medskip

Let $(\M,\E)$ be a $\tau$-rigid pair of $\mod\,\T$.  For
any $M\in\M$, take a projective
presentation
\begin{equation}\label{q1}
\begin{array}{l}
P_1\xrightarrow{\pi^M} P_0
\xrightarrow{} M
\xrightarrow{} 0
\end{array}
\end{equation}
such that the class $\{\, \pi^M \,|\, M \in \M \,\}$ has Property
(S). By Theorem \ref{b1}(i), there is a unique morphism
$f_M\colon R_1 \rightarrow R_0$ in $\R$ such that $\mathbb{H}( f_M ) =
\pi^M$. Moreover, $\mathbb{H}( \mathrm{cone}(f_M) ) \cong M$. Since
the projective
presentation (\ref{q1}) has Property (S), it
follows from Lemma \ref{b2} that the category
\[
  \X_1 := \{\, \mathrm{cone}(f_M) \;|\; M \in \M \,\}
\]
is a $\R[1]$-rigid subcategory and $\X_1\in\RR$.
\smallskip

Let $\X:=\X_1\vee\E[1]$. Since $\E\in\R$, we have $\E[1]\subseteq\R[1]\subseteq\RR$. Then $\X\subseteq\RR$.

Now we show that $\X$ is a $\R[1]$-rigid subcategory of $\C$.
Let $E\in\E\subseteq\R$. Since $\R$ is rigid, we have $$[\R[1]](\mathrm{cone}(f_M)\oplus E[1],E[2])=0.$$
Applying the functor $\Hom_{\C}(E,-)$ to the triangle
$R_1\xrightarrow{~f_M~}R_0\to \mathrm{cone}(f_M)\to R_1[1]$,
we have the following exact sequence
\[\Hom_{\C}(E,R_1)\xrightarrow{~f_M\circ~}\Hom_{\C}(E,R_0)\to\Hom_{\C}(E,\mathrm{cone}(f_M))\to \Hom_{\C}(E,R_1[1])=0,\]
which is isomorphic to $$P_1(E)\xrightarrow{~\pi^M~}P_0(E)\to M(E)\to 0.$$
The condition $\M(\E)=0$ implies that $\Hom_{\C}(E,\mathrm{cone}(f_M))=0$
and then $$[\R[1]](E[1],\mathrm{cone}(f_M)[1])=0.$$
Thus the assertion holds.
\smallskip

Now we show that $\Phi( \X ) =( \M , \E )$.
\smallskip

It is straightforward to check that $\R\cap\X_1[-1]=0$.
For any object $X\in\R\cap\X[-1]$, we can write $X=X_1[-1]\oplus E\in\R$, where $X_1\in\X_1$ and $E\in\E$.
Since $X_1[-1]\in\R\cap\X_1[-1]=0$, we have $X=E\in\E$. Thus we have $\R\cap\X[-1]\subseteq\E$.
By the definition of $\tau$-rigid pair, we have $\E\subseteq\R$.
Note that $\E\subseteq\X_1[-1]\vee\E=\X[-1]$, it follows that $\E\subseteq\R\cap\X[-1]$.
Hence $\R\cap\X[-1]=\E$.

It remains to show that $\mathbb{H}(\X)=\M$. Indeed, since $\E\subseteq\T$, we have
$$\mathbb{H}(\X)=\Hom_{\C}(\R,\X)=\Hom_{\C}(\R,\X_1)=\mathbb{H}(\X_1)=\M.$$

\textbf{Step 3:} The map $\Phi$ is injective.
\medskip

Let $\X$ and $\X'$ be two $\R[1]$-rigid subcategories of $\C$ which are contained in $\RR$ such that $\Phi(\X)=\Phi(\X')$. Let
$\X_1$ and $\X'_1$ be respectively the full subcategories of $\X$
and $\X'$ consisting of objects without direct summands in
$\R[1]$. Then $\X=\X_1\vee (\X\cap\R[1])$ and
$\X'=\X'_1\vee (\X'\cap\R[1])$. Since $\Phi(\X)=\Phi(\X')$,
it follows that $\mathbb{H}(\X_1)=\mathbb{H}(\X'_1)$ and
$\X\cap\R[1]=\X'\cap\R[1]$.
\smallskip

For any object $X_1\in\X_1$, there exists $X'_1\in\X'_1$ such that $\mathbb{H}(X_1)=\mathbb{H}(X'_1)$. By Theorem \ref{b1}(ii), there exists an isomorphism
$X_1\oplus Y[1]\simeq X'_1\oplus Z[1]$ for some $Y,Z\in\R$. Since $\C$ is Krull-Schmidt, we have
$X_1\simeq X'_1$. This implies that $\X_1\subseteq \X'_1$. Similarly, we obtain $\X'_1\subseteq \X_1$ and then $\X_1\simeq\X'_1$.
 Therefore $\X=\X'$. This shows that $\Phi$ is injective. \qed

\begin{thm}\label{b4}
Let $\C$ be a triangulated category and $\R$ a rigid category of $\C$.
The functor $$\mathbb{H}\colon \C\to \Mod\, \R$$ induces a bijection $$\Phi\colon \X\longmapsto \big(\mathbb{H}(\X),\R\cap\X[-1]\big)$$
from the first of the following sets to the second:
\begin{itemize}
\item[\emph{(I)}] two-term weak $\R[1]$-cluster tilting subcategories of $\C.$

\item[\emph{(II)}] Support $\tau$-tilting pairs of $\mod\, \R$.
\end{itemize}
\end{thm}

\proof  \textbf{Step 1:} The map $\Phi$ has values in support $\tau$-tilting pairs of $\mod\, \R$.
\smallskip

Assume that $\X$ is a weak $\R[1]$-cluster tilting subcategory of $\C$ which is contained in $\RR$. By Theorem \ref{b3}, we know that $\Phi(\X)$ is a $\tau$-rigid pair of $\mod\, \R$. Therefore $\R\cap\X[-1]\subseteq \Ker\  \mathbb{H}(\X).$

Let $R\in\R$ be an object of $\Ker \ \mathbb{H}(\X)$, i.e. $\Hom_{\C}( R,X )=0$ for each
$X\in\X$. This implies that $[\R[1]](X \oplus R[1], \X[1])=0$. Note that $[\R[1]]\big(\X, (X\oplus R[1])[1]\big)=0$. Since $\X$ is a $\R[1]$-cluster tilting, we have $X\oplus R[1]\in\X$ and then $R\subseteq\X[-1]$. Therefore $R\in \R\cap \X[-1]$.
This shows that $\Ker\  \mathbb{H}(\X) \subseteq \R\cap\X[-1].$ Hence
$$\Ker\ \mathbb{H}(\X) = \R\cap \X[-1].$$
By the definition of weak $\T[1]$-cluster tilting, for any $R\in\R$, there exists a triangle
$$R\xrightarrow{~f~}X_1\xrightarrow{~g~}X_2\xrightarrow{~h~}R[1],$$
where $X_1,X_2\in\X$.
Applying the functor $\mathbb{H}$ to the above triangle, we obtain an exact sequence
$$
  \mathbb{H}(R)
  \xrightarrow{~\mathbb{H}(f)~}\mathbb{H}(X_1)
  \longrightarrow \mathbb{H}(X_2)
  \rightarrow 0.
$$
For any morphism $a\colon T\to X$, where $X \in \X$, since $\X$ is $\R[1]$-rigid, we have $a\circ h[-1]=0$.
Then there exists a morphism $b\colon X_1\to X$ such that $a=bf$. This shows that $\Hom_{\C}(f, X)$
is a surjective. Thus there exists the following commutative diagram.
\[
  \xymatrix{
    \Hom_{\C}(X_1,X)\ar[r]^{\Hom_{\C}(f,\ X)}\ar[d] & \Hom_{\C}(R,X)\ar[r] \ar[d] &0\\
    \Hom_{\mod\,\R}(\mathbb{H}(X_1),\mathbb{H}(X))\ar[r]^{\circ\mathbb{H}(f)} & \Hom_{\mod\,\R}(\mathbb{H}(R),\mathbb{H}(X))
           }
\]
Using Theorem \ref{b1}(i), the right vertical map is an isomorphism. It follows that $\circ\mathbb{H}(f)$ is surjective, that is, $\mathbb{H}(f)$ is a left
$\mathbb{H}(\X)$-approximation of $\mathbb{H}(R)$.  This shows that $\Phi( \X )$
is a support $\tau$-tilting pair of $\mod\, \R$.
\medskip

\textbf{Step 2:} The map $\Phi$ is surjective.
\medskip

Let $(\M,\E)$ be a support $\tau$-tilting pair of $\mod\,\R$ and let $\X$ be the preimage of $(\M,\E)$ under $\Phi$ constructed in Theorem \ref{b3}. Since $\mathbb{H}(\X)=\M$ is a support $\tau$-tilting subcategory, for any $R\in\R$, there exists an exact sequence
$$
  \mathbb{H}(R)
  \xrightarrow{~\alpha~}\mathbb{H}(X_3)
  \rightarrow \mathbb{H}(X_4)
  \rightarrow 0.
$$
such that $X_3,X_4\in\X$ and $\alpha$ is a left
$\mathbb{H}(\X)$-approximation of $\mathbb{H}(R)$.  By Yoneda's lemma, there exists a unique morphism
$\beta\colon R\rightarrow X_3$ such that $\mathbb{H}( \beta ) = \alpha$.  Complete it to a triangle
\begin{equation}\label{q5}
\begin{array}{l}
R\s{\beta}\longrightarrow X_3 \s{\gamma}\longrightarrow Y_R \s{\delta}\longrightarrow R[1].
\end{array}
\end{equation}
Since $X_3\in\X\subseteq\RR$, there exists a triangle
$$R'\s{~}\longrightarrow R''\s{~}\longrightarrow X \s{~}\longrightarrow R'[1]$$
where $R',R''\in\R$. By the octahedral axiom,
we have the following commutative diagram
\begin{center}
$\xymatrix@!@C=0.5cm@R=0.5cm{
               &R'\ar[d]\ar@{=}[r]  & R'\ar[d]  \\
R\ar@{=}[d]\ar[r]  & X_3\ar[d]\ar[r]     & Y_R\ar[r]\ar[d] & R[1]\ar@{=}[d]  \\
R\ar[r]^{h=0\;\;}           & R''[1]\ar[d]\ar[r]    & Q\ar[d]\ar[r]    & R[1] \\
               & R'[1]\ar@{=}[r]      & R'[1]
  }$
\end{center}
of triangles. Since $\R$ is a rigid subcategory, we have $h=0$. It follow that $Q\simeq R''[1]\oplus R'[1]\in\R[1]$
implies $Y_R\in\RR$.
\smallskip

Let $\widetilde{\X}:= \X\vee\add\{\, Y_R
\,|\, R \in \R \,\}$ be the additive closure of $\X$ and $\{\, Y_R
\,|\, R\in \R \,\}$.  We obtain  $\widetilde{\X}\in\RR$.

We claim that $\widetilde{\X}$ is a weak $\R[1]$-cluster tilting
subcategory of $\C$ such that $\Phi(
\widetilde{\X} ) = ( \M , \E )$.
\smallskip

It is clear that $\R\subseteq\widetilde{\X}[-1]\ast\widetilde{\X}$. \   It remains to show that
$$\widetilde{\X}=\{M\in\RR\ |\ [\R[1]](M, \widetilde{\X}[1])=0=[\R[1]](\widetilde{\X}, M[1])\}.$$
Applying the functor $\mathbb{H}$ to the triangle (\ref{q5}), we see that $\mathbb{H}(Y_R)$ and $\mathbb{H}(X_4)$
are isomorphic in $\mod\,\R$. For any object $X\in\X$, we consider the following
commutative diagram.
\[
  \xymatrix{
    \Hom_{\C}(X_3,X) \ar[r]^{\Hom_{\C}(\beta,\ X)} \ar[d]_{\mathbb{H}(-)}
      & \Hom_{\C}(R,X)\ar[d]^{\simeq}\\
    \Hom_{\mod\,\R}( \mathbb{H}(X_3) , \mathbb{H}(X) ) \ar[r]^{\circ\alpha}
      & \Hom_{\mod\,\R}(\mathbb{H}(R),\mathbb{H}(X))}
\]
By Theorem \ref{b1}(i), the map $\mathbb{H}(-)$ is surjective and the right vertical map is an isomorphism.
since $\alpha$ is a left
$\mathbb{H}(\X)$-approximation of $\mathbb{H}(R)$, $\circ\alpha$ is also surjective. Therefore $\Hom_{\C}(\beta, X)$ is surjective too.
\smallskip

For any morphism $a\in[\R[1]](Y_R,X[1])$, since $\X$ is $\R[1]$-rigid, we have
$a\gamma=0$. So there exists a morphism $b\colon R[1]\to X[1]$ such that $a=b\delta$.
$$\xymatrix@!@C=0.5cm@R=0.3cm{
 R\ar[r]^{\beta} & X_3\ar[r]^{\gamma} & Y_R\ar[r]^{\delta\;\;} \ar[d]^{a}& R[1]\ar@{.>}[dl]^{b} \\
 && X[1]&}$$
Since $\Hom_{\C}(\beta,\ X)$ is surjective, there exists a morphism $c\colon X_3\to X$ such that
$c\beta=b[-1]$ and then $b=c[1]\circ \beta[1]$. It follows that
$a=b\delta=c[1]\circ(\beta[1]\delta)=0$.
This shows that
\begin{equation}\label{equi:mainthm1}
\begin{array}{l}
[\R[1]](Y_R,\X[1])=0.
\end{array}
\end{equation}

For any morphism $x\in[\R](X[-1],Y_R)$, we know that there exist two morphisms $x_1\colon X[-1]\to R_1$ and
$x_2\colon R_1\to Y_T$ such that $x=x_2x_1$, where $R_1\in\R$. Since $\R$ is rigid, we have $\delta x_2=0$. So there
exists a morphism $y\colon R_1\to X_3$ such that $x_2=\gamma y$.
$$\xymatrix@!@C=0.5cm@R=0.3cm{
&&X[-1]\ar[d]^{x_1}\\
&&R_1\ar[d]^{x_2}\ar@{.>}[dl]_{y}\\
R\ar[r]^{\beta} & X_3\ar[r]^{\gamma} & Y_R\ar[r]^{\delta\;\;}& R[1]
 }$$
Since $\X$ is $\R[1]$-rigid, we have
$x=x_2x_1=\gamma(yx_1)=0$. This shows that
\begin{equation}\label{equi:mainthm2}
\begin{array}{l}
[\R[1]](\X,Y_R[1])=0.
\end{array}
\end{equation}

For any $R'''\in\R$ and morphism $u\in [\R](Y_{R'''}[-1],Y_R)$, we know that there exist two morphisms $u_1\colon Y_{R'}[-1]\to T_2$ and
$u_2\colon R_2\to Y_R$ such that $u=u_2u_1$, where $R_2\in\R$.
Since $\T$ is rigid, we have $\delta u_2=0$. So there
exists a morphism $v\colon R_2\to X_3$ such that $u_2=\gamma v$.
$$\xymatrix@!@C=0.5cm@R=0.3cm{
&&Y_{R'''}[-1]\ar[d]^{u_1}\\
&&R_2\ar[d]^{u_2}\ar@{.>}[dl]_{v}\\
R\ar[r]^{\beta} & X_3\ar[r]^{\gamma} & Y_R\ar[r]^{\delta\;\;}& R[1]
 }$$
Since $[\R[1]](Y_R,\X[1])=0$, we have $v u_1=0$.
It follows that $u=u_2u_1=\gamma vu_1=0$. This shows that
\begin{equation}\label{equi:mainthm3}
\begin{array}{l}
[\R[1]](Y_{R'''},Y_R[1])=0.
\end{array}
\end{equation}
Using equalities (\ref{equi:mainthm1}), (\ref{equi:mainthm2}) and (\ref{equi:mainthm3}), we know that $\widetilde{\X}$ is a $\R[1]$-rigid subcategory.
\medskip

Now we show that $\{M\in\RR\ |\ [\R[1]](M, \widetilde{\X}[1])=0=[\R[1]](\widetilde{\X}, M[1])\}\subseteq\widetilde{\X}.$
\medskip

For any object $M\in\RR$, assume that $[\R[1]](M, \widetilde{\X}[1])=0=[\R[1]](\widetilde{\X}, M[1])$. Then there exists a triangle
$$R_5\s{f}\longrightarrow R_6 \s{g}\longrightarrow M \s{h}\longrightarrow R_5[1],$$
where $R_5,R_6\in\R$. By the above discussion, for object $R_6\in\R$, there exists a triangle
$$R_6\s{u}\longrightarrow X_6 \s{v}\longrightarrow Y_{R_6} \s{w}\longrightarrow R_6[1],$$
where $X_6\in\X$, $Y_{R_6}\in\widetilde{\X}$ and $u$ is a left $\X$-approximation of $R_6$.
For object $R_5\in\R$, there exists a triangle
$$R_5\s{u'}\longrightarrow X_5\s{v'}\longrightarrow Y_{R_5} \s{w'}\longrightarrow R_5[1],$$
where $X_5\in\X$, $Y_{R_5}\in\widetilde{\X}$ and $u'$ is a left $\X$-approximation of $R_5$.
By the octahedral axiom, we have a commutative diagram
$$\xymatrix{
R_5\ar[r]^{f}\ar@{=}[d]&R_6\ar[r]^{g}\ar[d]^{u}&M\ar[r]^{h}\ar[d]^{a}&R_5[1]\ar@{=}[d]\\
R_5\ar[r]^{x=uf}&X_6\ar[r]^{y}\ar[d]^{v}&N\ar[r]^{z}\ar[d]^{b}&R_5[1]\\
&Y_{R_6}\ar@{=}[r]\ar[d]^{w}&Y_{R_6}\ar[d]^{c}\\
&R_6[1]\ar[r]^{g[1]}&M[1]}$$
of triangles in $\C$.  We claim that $x$ is a left $\X$-approximation of $R_5$. Indeed, for any $d\colon R_5\to X$,
since $dh[-1]\in[\R](M[-1], \widetilde{\X})=0$, there exists a morphism $e\colon R_6\to X$ such that $d=ef$, where $X\in\X$.
$$\xymatrix{
    M[-1]\ar[r]^{\;\;\, h[-1]}&R_5 \ar[r]^{f} \ar[d]^{d} & R_6\ar@{.>}[dl]^{e} \ar[r]^{g} &M\ar[r]^{h}&R_5[1]\\
    &X&}$$
Since $u$ is a left $\X$-approximation of $R_6$, there exists a morphism $k\colon X_6\to X$ such that $ku=e$. It follows that $d=ef=kuf=kx$, as required.
\smallskip

Since $x$ is a left $\X$-approximation of $R_5$, by Lemma 1.4.3 in \cite{Ne},  we have the following commutative diagram
$$\xymatrix{
    R_5 \ar[r]^{x} \ar@{=}[d] & X_6\ar[r]^{y}\ar[d]^{\lambda} &N\ar[r]^{z}\ar@{.>}[d]^{\varphi}&R_5[1]\ar@{=}[d]\\
    R_5\ar[r]^{u'}&X_5\ar[r]^{v'}&Y_{R_5}\ar[r]^{w'}&R_5[1],}$$
where the middle square is homotopy cartesian and the differential $\partial=x[1]\circ w'$, that is, there exists a triangle
$$X_6\xrightarrow{\binom{-y}{\lambda}} N\oplus X_5 \xrightarrow{(\varphi,\ v')} Y_{R_5} \xrightarrow{~~\partial~~}X_6[1].$$
Note that $\partial\in[\R[1]](\widetilde{\X},\widetilde{\X}[1])=0$.
Thus we have $N\oplus X_5\simeq X_6\oplus Y_{R_5}\in\widetilde{\X}$, which implies $N\in\widetilde{\X}$.
 Since $c=g[1]w\in [\T[1]](\widetilde{\X},M[1])=0$, we know that the triangle
$$\xymatrix{M\ar[r]^{a}&N\ar[r]^{b}&Y_{R_6}\ar[r]^c&M[1]}$$
splits.  Hence $M$ is a direct summand of $N$ and then $M\in\widetilde{\X}$.
\medskip

This shows that $\widetilde{\X}=\{M\in\C\ |\ [\R[1]](M, \widetilde{\X}[1])=0=[\R[1]](\widetilde{\X}, M[1])\}.$
\smallskip

For any object $R\in\R$, $\mathbb{H}(Y_R)\simeq \mathbb{H}(X_4)$. Therefore $$\mathbb{H}(\widetilde{\X})\simeq \mathbb{H}(\X)\simeq \M.$$
Since $\R\cap \widetilde{\X}[-1]\supseteq \R\cap \X[-1]=\E$ and $\R\cap \widetilde{\X}[-1]\subseteq
\Ker\ \mathbb{H}(\X)=\E$, we have $$\R\cap \widetilde{\X}[-1]=\E.$$
This shows that $\Phi$ is surjective.
\medskip

\textbf{Step 3:} The map $\Phi$ is injective.
\smallskip

This follows from  Step 3 in Theorem \ref{b3}.  \qed

\medskip
\begin{rem}
If $\R=\add R$, Theorem \ref{b3} and Theorem \ref{b4} are just the Theorem 2.5 in \emph{\cite{FGL}}.
\end{rem}

Now we recover the following result in \cite[Theorem 4.13]{YZZ}.

\begin{cor}
Let $\R$ be a cluster tilting subcategory of $\C$. Then there exists a bijection between the set of weak $\R[1]$-cluster tilting subcategories and the set of support $\tau$-tilting pairs of $\mod\R$.
\end{cor}

\section{Relative rigid subcategories and presilting subcategories}
\setcounter{equation}{0}

We write the vanishing condition $\Hom(X,Y[i])=0$ for any $i>m$ by $\Hom(X,Y[>m])=0$.

\begin{defn}\emph{\cite[Definition 2.1]{AI} and \cite[Definition 1.2]{IJY}}\label{d3}
Let $\C$ be a triangulated category.
\begin{itemize}
\item[\emph{(i)}] A subcategory $\X$ of $\C$ is called a \emph{presilting subcategory} if $\emph{\Hom}_{\C}(\X,\X[>0])=0$.

\item[\emph{(ii)}] A presilting subcategory $\X\subseteq\C$ is a \emph{silting subcategory} if $\thick(\X)=\C$.
Here $\thick(\X)$ denotes the smallest thick subcategory of $\C$ containing $\X$.

\item[\emph{(iii)}] We say that an object $X\in\C$ is \emph{presilting }(respectively, \emph{silting}) if so is $\add X$.
\end{itemize}
\end{defn}

Let $\C$ be a triangulated category and $\S$ a silting subcategory of $\C$. A silting subcategory $\X$ of $\C$ is called \emph{two-term} with respect to $\S$ if $\X\subseteq\S\ast\S[1]$, see \cite{IJY}.

\begin{lem}\emph{\cite[Corollary 2.4]{IJY}}\label{d4}
Let $\C$ be a triangulated category, $\S$ and $\S'$ be two silting subcategories of $\C$.
If $\S'$ is two-term with respect to $\S$, then $\S$ is two-term with respect to $\S'[-1]$.
\end{lem}

The proof of the following lemma is the same as Lemma 3.3 in \cite{FGL} by adapting the objects to subcategories. We omit the proof and refer to \cite{FGL}.

\begin{lem}\label{d5}
Let $\C$ be a triangulated category and $\R$ a presilting subcategory of $\C$.
If a subcategory $\X$ of $\C$ which is contained in $\RR$,  then $\emph{\Hom}_{\C}(\X,\X[>1])=0$.
\end{lem}

\begin{thm}\label{main2}
Let $\C$ be a triangulated category and $\X$ a subcategory of $\C$ which is contained in $\RR$.
\begin{enumerate}
\item[\emph{(1)}] If $\R$ is a presilting subcategory of $\C$,  then the following statements are equivalent:
\begin{itemize}
\item[\emph{(a)}] $\X$ is an  $\R[1]$-rigid subcategory;
\item[\emph{(b)}] $\X$ is a rigid subcategory;
\item[\emph{(c)}] $\X$ is a presilting subcategory.
\end{itemize}
\item[\emph{(2)}]  If $\R$ is a silting subcategory of $\C$,  then $\X$ is a weak $\R[1]$-cluster tilting subcategory of $\C$ if and only if $\X$ is a silting subcategory of $\C$.
 \end{enumerate}
\end{thm}

\proof (1)~ The proof is the same as Theorem 3.4 (1) in \cite{FGL}, we omit it.
\smallskip

(2) Now we assume that $\R$ is a silting subcategory of $\C$.
\smallskip

If $\X$ is a weak $\R[1]$-cluster tilting subcategory, we have that $\X$ is $\R[1]$-rigid. By (1), we know that
$\X$ is presilting. By the definition of $\R[1]$-cluster tilting, we have $\R\subseteq\X[-1]\ast\X$.
It follows that $\X$ generates $\R$. On the other hand, $\R$ generates $\C$ since $\R$ is a silting subcategory.
Hence $\X$ generates $\C$ and then $\X$ is silting.  \qed

Conversely, if $\X$ is a silting subcategory of $\C$, since $\R$ is silting and $\X\subseteq\RR$, we have
$\R\subseteq\X[-1]\ast\X$ by Lemma \ref{d4}.

Since $\X$ is silting, it is clear that $\X$ is $\R[1]$-rigid. In other words, we have
$$\X\subseteq \{M\in\RR\ |\ [\R[1]](\X, M[1])=0= [\R[1]](M, \X[1])\}.$$

For any object $M\in \{M\in\RR\ |\ [\R[1]](\X, M[1])=0= [\R[1]](M, \X[1])\}$,
there exists a triangle
$$R_1\xrightarrow{~f~}R_0\xrightarrow{~g~}X\xrightarrow{~h~}R_1[1],$$
where $R_0,R_1\in\R$. Since $R_0\in\R\subseteq\X[-1]\ast\X$, there exists a triangle
$$R_0\xrightarrow{~u~}X_1\xrightarrow{~v~}X_2\xrightarrow{~w~}R_0[1],$$
where $X_1,X_2\in\X$. Since $\X$ is $\R[1]$-rigid, we obtain that $u$ is a left $\X$-approximation of $R_0$.
By the octahedral axiom, we have a commutative diagram
$$\xymatrix{
R_1\ar[r]^{f}\ar@{=}[d]&R_0\ar[r]^{g}\ar[d]^{u}&M\ar[r]^{h}\ar[d]^{a}&R_1[1]\ar@{=}[d]\\
R_1\ar[r]^{x=uf}&X_1\ar[r]^{y}\ar[d]^{v}&N\ar[r]^{z}\ar[d]^{b}&R_1[1]\\
&X_2\ar@{=}[r]\ar[d]^{w}&X_2\ar[d]^{c}\\
&R_0[1]\ar[r]^{g[1]}&M[1]}$$
of triangles. We claim that $x$ is a left $\X$-approximation of $R_1$.
Indeed, for any
morphism $\alpha\colon R_1\to X'$, where $X'\in\X$, since $\alpha\circ h[-1]\in[\R](M[-1],\X)=0$, there exists a morphism $\beta\colon R_0\to X'$ such that $\alpha=\beta f$. Since $u$ is a left $\X$-approximation of $R_0$ and $X'\in\X$,  there exists a morphism $\gamma\colon X_1\to X'$ such that
$\beta=\gamma u$ and then $\alpha=\gamma(uf)=\gamma x$. This shows that $x$ is a left $\X$-approximation of $R_1$.

Since $R_1\in\R\subseteq\X[-1]\ast\X$, there exists a triangle
$$R_1\xrightarrow{~u'~}X_3\xrightarrow{~v'~}X_4\xrightarrow{~w'~}R_1[1],$$
where $X_3,X_4\in\X$. By \cite[Lemma 1.4.3]{Ne} and $x$ is a left $\X$-approximation of $R_1$,  we have the following commutative diagram
$$\xymatrix{
    R_1 \ar[r]^{x} \ar@{=}[d] & X_1\ar[r]^{y}\ar@{.>}[d]^{\lambda} &N\ar[r]^{z\;\,}\ar@{.>}[d]^{\varphi}&R_1[1]\ar@{=}[d]\\
    R_1\ar[r]^{u'}&X_3\ar[r]^{v'}&X_4\ar[r]^{w'\;\,}&R_1[1],}$$
where the middle square is homotopy cartesian and the differential $\partial=x[1]\circ w'$, that is, there exists a triangle
$$X_1\xrightarrow{\binom{-y}{\lambda}} N\oplus X_3 \xrightarrow{(\varphi,\ v')} X_4 \xrightarrow{~~\partial~~}X_1[1].$$
Note that $\partial\in[\R[1]](\X,M[1])=0$.
Thus we have $N\oplus X_3\simeq X_1\oplus X_4\in\X$, which implies $N\in\X$.
 Since $c=g[1]\circ w\in [\R[1]](\X,M[1])=0$, we know that the triangle
$$\xymatrix{M\ar[r]^{a}&N\ar[r]^{b}&Y_{T_6}\ar[r]^c&M[1]}$$
splits.  Hence $M$ is a direct summand of $N$ and then $M\in\X$.

This shows that $\X$ is weak $\R[1]$-cluster tilting.  \qed
\medskip

Combining Theorem \ref{b3}, Theorem \ref{b4} with Theorem \ref{main2}, we obtain the following bijection, which is due to Iyama-J{\o}rgensen-Yang, see \cite[Theorem 3.3]{IJY}.
\begin{cor}
Let $\R$ be a silting subcategory of $\C$. There exists a bijection between the set of presilting subcategories of $\C$ which are contained in $\RR$ and the set of $\tau$-rigid pairs of $\mod\R$, which induces a one-to-one correspondence between the set of  silting subcategories which are contained in $\RR$ and the set of support $\tau$-tilting pairs of $\mod\R$.
\end{cor}

Panyue Zhou\\
College of Mathematics, Hunan Institute of Science and Technology, 414006 Yueyang, Hunan, People's Republic of China.\\[1mm]
Email: panyuezhou@163.com\\[3mm]
Bin Zhu\\
Department of Mathematical Sciences,
Tsinghua University,
100084, Beijing,
People's Republic of China\\[1mm]
Email: bzhu@math.tsinghua.edu.cn


\begin{thebibliography}{99}
\bibitem[Au]{Au}
M.\ Auslander.
\newblock Representation theory of Artin
algebras I.
\newblock Comm.\ Algebra, 1, 177-268, 1974.


\bibitem[AI]{AI} T. Aihara and O. Iyama. Silting mutation in triangulated categories. J. Lond. Math. Soc. 85(3), 633-668, 2012.

\bibitem[AV]{AV} J. Stovicek and A. Van Roosmalen. \newblock $2$-Calabi-Yau categories with a directed cluster-tilting subcategory.
\newblock arXiv: 1611. 03836v1, 2016.

\bibitem[AIR]{AIR}
T. Adachi, O. Iyama and I.  Reiten. $\tau$-tilting theory. Compos. Math. 150(3), 415-452, 2014.

\bibitem[BMRRT]{BMRRT} A. Buan, R. Marsh, M. Reineke, I. Reiten and G. Todorov. \newblock Tilting theory and cluster combinatorics.  \newblock Adv. Math. 204(2), 572-618, 2006.


\bibitem[CZZ] {CZZ}
W. Chang, J. Zhang and B. Zhu. On support $\tau$-tilting modules over endomorphism algebras of rigid objects,
Acta Math. Sin. (Engl. Ser.), 31(9), 1508-1516, 2015.



\bibitem[FGL]{FGL} C. Fu, S. Geng and P. Liu. Relative rigid objects in triangulated categories.
\newblock arXiv: 1808.04297, 2018.

\bibitem[IY]{IY}
O. Iyama and Y. Yoshino. Mutation in triangulated categories and rigid Cohen-Macaulay modules. Invent. Math. 172, 117-168, 2008.


\bibitem[IJY]{IJY}
O. Iyama, P. J{\o}rgensen and D. Yang. Intermediate co-t-structures, two-term silting objects, $\tau$-tilting modules, and torsion classes,  Algebra and Number Theory, 8(10), 2413-2431, 2014.


\bibitem[KR]{KR}
B. Keller and I. Reiten.
\newblock Cluster-tilted algebras are Gorenstein and stably Calabi-Yau.
\newblock Adv. Math. 211, 123-151, 2007.

\bibitem[KZ]{KZ}
S. Koenig and B. Zhu.
\newblock From triangulated categories to abelian categories: cluster tilting in a general framework.
\newblock Math. Z. 258, 143-160, 2008.



\bibitem[LX]{LX}
 P. Liu and Y. Xie. On the relation between maximal rigid objects and $\tau$-tilting modules.
Colloq. Math. 142(2), 169-178, 2016.


\bibitem[Ne]{Ne} A.\ Neeman. \newblock Triangulated Categories. Volume 148. Princeton University Press, 2001.


\bibitem[YZ]{YZ}
W. Yang and B. Zhu. \newblock Relaive cluster tilting objects in triangulated categories. arXiv:1504.00093,  2015,  to
appear in Trans. Amer. Math. Soc.

\bibitem[YZZ]{YZZ}
W. Yang, P. Zhou and  B. Zhu. \newblock Triangulated categories with cluster-tilting subcategories,  arXiv:1711.04290, 2017.

\bibitem[ZZ]{ZZ} Y.\ Zhou and B.\ Zhu. \newblock Maximal rigid subcategories in $2$-Calabi-Yau triangulated categories.
\newblock J. Algebra, 348: 49-60, 2011.



\end{thebibliography}
\end{document}